\documentclass{article}
\usepackage{amssymb}
\usepackage{amsmath}
\usepackage{amstext}
\usepackage{graphicx}

\newtheorem{theorem}{Theorem}

\newtheorem{axiom}[theorem]{Axiom}

\newtheorem{corollary}[theorem]{Corollary}

\newtheorem{definition}[theorem]{Definition}

\newtheorem{proposition}[theorem]{Proposition}
\newtheorem{remark}[theorem]{Remark}

\input tcilatex
\begin{document}

\title{Basic properties of ultrafunctions}
\author{Vieri Benci\thanks{
Dipartimento di Matematica, Universit\`{a} degli Studi di Pisa, Via F.
Buonarroti 1/c, 56127 Pisa, ITALY and Department of Mathematics, College of
Science, King Saud University, Riyadh, 11451, SAUDI ARABIA. e-mail: \texttt{%
benci@dma.unipi.it}} \and Lorenzo Luperi Baglini\thanks{%
Dipartimento di Matematica, Universit\`{a} degli Studi di Pisa, Via F.
Buonarroti 1/c, 56127 Pisa, ITALY, e-mail:\texttt{\
lorenzo.luperi@for.unipi.it}}}
\maketitle

\begin{center}

\emph{Dedicated to Bernard Ruf in occasion of his 60th birthday.}

\end{center}

\medskip

\begin{abstract}
Ultrafunctions are a particular class of functions defined on a
non-Archimedean field $\mathbb{R}^{\ast }\supset \mathbb{R}$. They provide
generalized solutions to functional equations which do not have any
solutions among the real functions or the distributions. In this paper we
analyze sistematically some basic properties of the spaces of ultrafunctions.

\medskip 


\medskip

\end{abstract}

\tableofcontents

\section{Introduction}

In some recent papers the notion of ultrafunction has been introduced (\cite%
{ultra}, \cite{belu2012}). Ultrafunctions are a particular class of
functions defined on a non-Archimedean field $\mathbb{R}^{\ast }\supset 
\mathbb{R}$. We recall that a non-Archimedean field is an ordered field
which contain infinite and infinitesimal numbers.

To any continuous function $f:\mathbb{R}^{N}\rightarrow \mathbb{R}$ we
associate in a canonical way an ultrafunction $\widetilde{f}:\left( \mathbb{R%
}^{\ast }\right) ^{N}\rightarrow \mathbb{R}^{\ast }$ which extends $f;$ more
exactly, to any functional vector space $V(\Omega )\subseteq L^{2}(\Omega
)\cap \mathcal{C}(\overline{\Omega }),$ we associate a space of
ultrafunctions $\widetilde{V}(\Omega ).\ $The ultrafunctions are much more
than the functions and among them we can find solutions of functional
equations which do not have any solutions among the real functions or the
distributions.

A typical example of this situation is analyzed in \cite{belu2012} where a
simple Physical model is studied. In this problem there is a material point
interacting with a field and, as it usually happens, the energy is infinite.
Therefore the need to use infinite numbers arises naturally. Other
situations in which infinite and infinitesimal numbers appear in a natural
way are studied in \cite{BGG} and in \cite{BHW}.

\bigskip

In this paper we analyze systematically some basic properties of the spaces
of ultrafunctions $\widetilde{V}(\Omega )$. In particular we will show that:

\begin{itemize}
\item to any measurable function $f$ we can associate an unique
ultrafunction $\widetilde{f}$ such that $f(x)=\widetilde{f}(x)$ if $f$ is
continuous in a neighborhood of $x;$

\item to every distribution $T$ we can associate an ultrafunction $%
\widetilde{T}(x)$ such that $\forall \varphi \in \mathcal{D},$ $\left\langle
T,\varphi \right\rangle =\int^{\ast }\widetilde{T}(x)\widetilde{\varphi }%
(x)dx$ where $\int^{\ast }$ is a suitable extension of the integral to the
ultrafunctions;

\item the vector space of ultrafunctions $\widetilde{V}(\Omega )$ is
hyperfinite, namely it shares many properties of finite vector spaces (see
Sec. \ref{HE});

\item the vector space of ultafunctions $\widetilde{V}(\Omega )$ has an
hyperfinite basis $\left\{ \delta _{a}(x)\right\} _{a\in \Sigma }$ where $%
\delta _{a}$ is the "Dirac ultrafunction in $a$" (see Def. \ref{dede}) and $%
\Sigma \subset \left( \mathbb{R}^{\ast }\right) ^{N}$ is a suitable set;

\item any ultrafunction $u$ can be represented as follows:%
\begin{equation*}
u(x)=\sum_{q\in \Sigma }u(q)\sigma _{q}(x),
\end{equation*}%
where $\left\{ \sigma _{a}(x)\right\} _{a\in \Sigma }$ is the dual basis of $%
\left\{ \delta _{a}(x)\right\} _{a\in \Sigma };$

\item any operator $F:V\left( \Omega \right) \rightarrow \mathcal{D}^{\prime
}\left( \Omega \right) ,$ can be extended to an operator%
\begin{equation*}
\widetilde{F}:\widetilde{V}\left( \Omega \right) \rightarrow \widetilde{V}%
\left( \Omega \right) ;
\end{equation*}%
the extension of the derivative and the Fourier transform will be analyzed
in some detail.
\end{itemize}

\bigskip

The techniques on which the notion of ultrafunction is based are related to
Non Archimedean Mathematics (NAM) and to Nonstandard Analysis (NSA). The
first section of this paper is devoted to a relatively elementary
presentation of the basic notions of NAM and NSA inspired by \cite{BDN2003}
and \cite{benci99}. Some technicalities have been avoided by presenting the
matter in an axiomatic way. Of course, it is necessary to prove the
consistency of the axioms. This is done in the appendix; however in the
appendix we have assumed the reader to be familiar with NSA.

\subsection{Notations\label{not}}

Let $\Omega $\ be a subset of $\mathbb{R}^{N}$: then

\begin{itemize}
\item $\mathcal{C}\left( \Omega \right) $ denotes the set of continuous
functions defined on $\Omega \subset \mathbb{R}^{N};$

\item $\mathcal{C}_{0}\left( \Omega \right) $ denotes the set of continuous
functions in $\mathcal{C}\left( \Omega \right) $ having compact support in $%
\Omega ;$

\item $\mathcal{C}^{k}\left( \Omega \right) $ denotes the set of functions
defined on $\Omega \subset \mathbb{R}^{N}$ which have continuous derivatives
up to the order $k;$

\item $\mathcal{D}\left( \Omega \right) $ denotes the set of the infinitely
differentiable functions with compact support defined on $\Omega \subset 
\mathbb{R}^{N};\ \mathcal{D}^{\prime }\left( \Omega \right) $ denotes the
topological dual of $\mathcal{D}\left( \Omega \right) $, namely the set of
distributions on $\Omega ;$

\item $H^{1,p}(\Omega )$ is the usual Sobolev space defined as the set of
functions in $L^{p}\left( \Omega \right) $ such that $\nabla u\in
L^{p}\left( \Omega \right) ^{N};$

\item $H^{1}(\Omega )=H^{1,2}(\Omega )$

\item for any $\xi \in \left( \mathbb{R}^{N}\right) ^{\ast },\rho \in 
\mathbb{R}^{\ast }$, we set $\mathfrak{B}_{\rho }(\xi )=\left\{ x\in \left( 
\mathbb{R}^{N}\right) ^{\ast }:\ |x-\xi |<\rho \right\} $;

\item $\mathfrak{supp}(f)=\overline{\left\{ x\in \mathbb{R}^{N}:f(x)\neq
0\right\} };$

\item $\mathfrak{mon}(x)=\{y\in \mathbb{R}^{N}:x\sim y\};$

\item $\mathfrak{gal}(x)=\{y\in \mathbb{R}^{N}:x\sim _{f}y\}.$
\end{itemize}

\section{$\Lambda $-theory\label{lt}}

In this section we present the basic notions of Non Archimedean Mathematics
and of Nonstandard Analysis following a method inspired by \cite{BDN2003}
(see also \cite{ultra} and \cite{belu2012}).

\subsection{Non Archimedean Fields\label{naf}}

Here, we recall the basic definitions and facts regarding non-Archimedean
fields. In the following, ${\mathbb{K}}$ will denote an ordered field. We
recall that such a field contains (a copy of) the rational numbers. Its
elements will be called numbers.

\begin{definition}
Let $\mathbb{K}$ be an ordered field. Let $\xi \in \mathbb{K}$. We say that:

\begin{itemize}
\item $\xi $ is infinitesimal if, for all positive $n\in \mathbb{N}$, $|\xi
|<\frac{1}{n}$;

\item $\xi $ is finite if there exists $n\in \mathbb{N}$ such as $|\xi |<n$;

\item $\xi $ is infinite if, for all $n\in \mathbb{N}$, $|\xi |>n$
(equivalently, if $\xi $ is not finite).
\end{itemize}
\end{definition}

\begin{definition}
An ordered field $\mathbb{K}$ is called Non-Archimedean if it contains an
infinitesimal $\xi \neq 0$.
\end{definition}

It's easily seen that all infinitesimal are finite, that the inverse of an
infinite number is a nonzero infinitesimal number, and that the inverse of a
nonzero infinitesimal number is infinite.

\begin{definition}
A superreal field is an ordered field $\mathbb{K}$ that properly extends $%
\mathbb{R}$.
\end{definition}

It is easy to show, due to the completeness of $\mathbb{R}$, that there are
nonzero infinitesimal numbers and infinite numbers in any superreal field.
Infinitesimal numbers can be used to formalize a new notion of "closeness":

\begin{definition}
\label{def infinite closeness} We say that two numbers $\xi, \zeta \in {%
\mathbb{K}}$ are infinitely close if $\xi -\zeta $ is infinitesimal. In this
case, we write $\xi \sim \zeta $.
\end{definition}

Clearly, the relation "$\sim $" of infinite closeness is an equivalence
relation.

\begin{theorem}
If $\mathbb{K}$ is a superreal field, every finite number $\xi \in \mathbb{K}
$ is infinitely close to a unique real number $r\sim \xi $, called the 
\textbf{shadow} or the \textbf{standard part} of $\xi $.
\end{theorem}

Given a finite number $\xi $, we denote its shadow as $sh(\xi )$, and we put 
$sh(\xi )=+\infty $ ($sh(\xi )=-\infty $) if $\xi \in \mathbb{K}$ is a
positive (negative) infinite number.\newline

\begin{definition}
Let $\mathbb{K}$ be a superreal field, and $\xi \in \mathbb{K}$ a number.
The \label{def monad} monad of $\xi $ is the set of all numbers that are
infinitely close to it:%
\begin{equation*}
\mathfrak{m}\mathfrak{o}\mathfrak{n}(\xi )=\{\zeta \in \mathbb{K}:\xi \sim
\zeta \},
\end{equation*}%
and the galaxy of $\xi $ is the set of all numbers that are finitely close
to it: 
\begin{equation*}
\mathfrak{gal}(\xi )=\{\zeta \in \mathbb{K}:\xi -\zeta \ \text{is\ finite}\}
\end{equation*}
\end{definition}

By definition, it follows that the set of infinitesimal numbers is $%
\mathfrak{mon}(0)$ and that the set of finite numbers is $\mathfrak{gal}(0)$.

\subsection{The $\Lambda $-limit\label{OL}}

In this section we will introduce a superreal field $\mathbb{K}$ and we will
analyze its main properties by mean of the $\Lambda $-theory (see also \cite%
{ultra}, \cite{belu2012}).

\bigskip

$\mathbb{U}$ will denote our "mathematical universe". For our applications a
good choice of $\mathbb{U}$ is given by the superstructure on $\mathbb{R}$:

\begin{equation*}
\mathbb{U}=\dbigcup_{n=0}^{\infty }\mathbb{U}_{n}
\end{equation*}%
where $\mathbb{U}_{n}$ is defined by induction as follows:%
\begin{eqnarray*}
\mathbb{U}_{0} &=&\mathbb{R}\text{;} \\
\mathbb{U}_{n+1} &=&\mathbb{U}_{n}\cup \mathcal{P}\left( \mathbb{U}%
_{n}\right) .
\end{eqnarray*}%
Here $\mathcal{P}\left( E\right) $ denotes the power set of $E.$ Identifying
the couples with the Kuratowski pairs and the functions and the relations
with their graphs, it follows that{\ }$\mathbb{U}$ contains almost every
usual mathematical object. Given the universe $\mathbb{U}$, we denote by $%
\mathcal{F}$ the family of finite subsets of $\mathbb{U}.$ Clearly $\left( 
\mathcal{F},\subset \right) $ is a directed set and, as usual, a function $%
\varphi :\mathcal{F}\rightarrow E$ will be called \textit{net }(with values
in $E$).\newline
We present axiomatically the notion of $\Lambda $-limit:

\bigskip {\Large Axioms of\ the }$\Lambda ${\Large -limit}

\begin{itemize}
\item \textsf{(}$\Lambda $-\textsf{1)}\ \textbf{Existence Axiom.}\ \textit{%
There is a superreal field} $\mathbb{K}\supset \mathbb{R}$ \textit{such that
every net }$\varphi :\mathcal{F}\rightarrow \mathbb{R}$\textit{\ has a
unique limit }$L\in \mathbb{K}{\ }($\textit{called the} "$\Lambda $-limit" 
\textit{of}\emph{\ }$\varphi .)$ \textit{The} $\Lambda $-\textit{limit of $%
\varphi $ will be denoted as} 
\begin{equation*}
L=\lim_{\lambda \uparrow {\mathbb{U}}}\varphi (\lambda ).
\end{equation*}%
\textit{Moreover we assume that every}\emph{\ }$\xi \in \mathbb{K}$\textit{\
is the }$\Lambda $-\textit{limit\ of some real function}\emph{\ }$\varphi :%
\mathcal{F}\rightarrow \mathbb{R}$\emph{. }

\item ($\Lambda $-2)\ \textbf{Real numbers axiom}. \textit{If }$\varphi
(\lambda )$\textit{\ is} \textit{eventually} \textit{constant}, \textit{%
namely} $\exists \lambda _{0}\in \mathcal{F},r\in \mathbb{R}$ such that $%
\forall \lambda \supset \lambda _{0},\ \varphi (\lambda )=r,$ \textit{then}%
\begin{equation*}
\lim_{\lambda \uparrow {\mathbb{U}}}\varphi (\lambda )=r.
\end{equation*}

\item ($\Lambda $-3)\ \textbf{Sum and product Axiom}.\ \textit{For all }$%
\varphi ,\psi :\mathcal{F}\rightarrow \mathbb{R}$\emph{: }%
\begin{eqnarray*}
\lim_{\lambda \uparrow {\mathbb{U}}}\varphi (\lambda )+\lim_{\lambda
\uparrow {\mathbb{U}}}\psi (\lambda ) &=&\lim_{\lambda \uparrow {\mathbb{U}}%
}\left( \varphi (\lambda )+\psi (\lambda )\right) ; \\
\lim_{\lambda \uparrow {\mathbb{U}}}\varphi (\lambda )\cdot \lim_{\lambda
\uparrow {\mathbb{U}}}\psi (\lambda ) &=&\lim_{\lambda \uparrow {\mathbb{U}}%
}\left( \varphi (\lambda )\cdot \psi (\lambda )\right) .
\end{eqnarray*}
\end{itemize}

\begin{theorem}
\label{brufolo}The set of axioms $\{$($\Lambda $-1)\textsf{,}($\Lambda $-2),(%
$\Lambda $-3)$\}$ is consistent.
\end{theorem}

Theorem \ref{brufolo} will be proved in the Appendix.

\bigskip

Now we want to define the $\Lambda $-limit of any bounded net of
mathematical objects in {$\mathbb{U}$} (a net $\varphi :\mathcal{F}%
\rightarrow ${$\mathbb{U}$} is called bounded if there exists $n$ such that $%
\forall \lambda \in \mathcal{F},\varphi (\lambda )\in ${$\mathbb{U}$}$_{n}$%
). To this aim, consider a net%
\begin{equation}
\varphi :\mathcal{F}\rightarrow {\mathbb{U}}_{n}.  \label{net}
\end{equation}%
We will define $\lim\limits_{\lambda \uparrow {\mathbb{U}}}\varphi (\lambda
) $ by induction on $n$. For $n=0,$ $\lim\limits_{\lambda \uparrow {\mathbb{U%
}}}\varphi (\lambda )$ is defined by the axioms \textsf{(}$\Lambda $-\textsf{%
1),}($\Lambda $-2),($\Lambda $-3); so by induction we may assume that the
limit is defined for $n-1$ and we define it for the net (\ref{net}) as
follows:%
\begin{equation*}
\lim_{\lambda \uparrow {\mathbb{U}}}\varphi (\lambda )=\left\{ \lim_{\lambda
\uparrow {\mathbb{U}}}\psi (\lambda )\ |\ \psi :\mathcal{F}\rightarrow 
\mathcal{\mathbb{U}}_{n-1}\text{ and}\ \forall \lambda \in \mathcal{F},\
\psi (\lambda )\in \varphi (\lambda )\right\} .
\end{equation*}

\begin{definition}
A mathematical entity (number, set, function or relation) which is the $%
\Lambda $-limit of a net is called \textbf{internal}.
\end{definition}

\subsection{Natural extensions of sets and functions\label{qualified}}

\begin{definition}
The \textbf{natural extension }of a set $E\subset \mathbb{R}$ is given by%
\begin{equation*}
E^{\ast }:=\lim_{\lambda \uparrow {\mathbb{U}}}c_{E}(\lambda )=\ \left\{
\lim_{\lambda \uparrow {\mathbb{U}}}\psi (\lambda )\ |\ \psi (\lambda )\in
E\right\}
\end{equation*}%
where $c_{E}(\lambda )$ is the net identically equal to $E$.
\end{definition}

This definition, combined with axiom ($\Lambda $-1$)$, entails that 
\begin{equation*}
\mathbb{K}=\mathbb{R}^{\ast }.
\end{equation*}

In this context a function $f$ can be identified with its graph; then the
natural extension of a function is well defined. Moreover we have the
following result:

\begin{theorem}
The \textbf{natural extension} of a function%
\begin{equation*}
f:E\rightarrow F
\end{equation*}%
is a function 
\begin{equation*}
f^{\ast }:E^{\ast }\rightarrow F^{\ast }
\end{equation*}%
and for every net $\varphi :\mathcal{F}\cap \mathcal{P}\left( E\right)
\rightarrow E,$ and every function $f:E\rightarrow F$, we have that%
\begin{equation*}
\lim_{\lambda \uparrow {\mathbb{U}}}\ f(\varphi (\lambda ))=f^{\ast }\left(
\lim_{\lambda \uparrow {\mathbb{U}}}\varphi (\lambda )\right) .
\end{equation*}
\end{theorem}

When dealing with functions, sometimes the "$\ast $" will be omitted if the
domain of the function is clear from the context. For example, if $\eta \in 
\mathbb{R}^{\ast }$ is an infinitesimal, then clearly $e^{\eta }$ denotes $%
\exp ^{\ast }(\eta ).$

The following theorem is a fundamental tool in using the $\Lambda $-limit:

\begin{theorem}
\label{limit}\textbf{(Leibnitz Principle)} Let $\mathcal{R}$ be a relation
in {$\mathbb{U}$}$_{n}$ for some $n\geq 0$ and let $\varphi $,$\psi :%
\mathcal{F}\rightarrow {\mathbb{U}}_{n}$. If 
\begin{equation*}
\forall \lambda \in \mathcal{F},\ \varphi (\lambda )\mathcal{R}\psi (\lambda
)
\end{equation*}%
then%
\begin{equation*}
\left( \underset{\lambda \uparrow {\mathbb{U}}}{\lim }\varphi (\lambda
)\right) \mathcal{R}^{\ast }\left( \underset{\lambda \uparrow {\mathbb{U}}}{%
\lim }\psi (\lambda )\right) .
\end{equation*}
\end{theorem}

When $\mathcal{R}$ is $\in $ or $\mathcal{=}$ we will not use the symbol $%
\ast $ to denote their extensions, since their meaning is unaltered in $%
\mathbb{R}^{\ast }.$

\subsection{Hyperfinite extensions\label{HE}}

\begin{definition}
An internal set is called \textbf{hyperfinite} if it is the $\Lambda $-limit
of a net $\varphi :\mathcal{F}\rightarrow \mathcal{F}$.
\end{definition}

\begin{definition}
Given any set $E\in $ {$\mathbb{U}$}, the hyperfinite extension of $E$ is
defined as follows:%
\begin{equation*}
E^{\circ }:=\ \underset{\lambda \uparrow {\mathbb{U}}}{\lim }(E\cap \lambda
).
\end{equation*}
\end{definition}

All the internal finite sets are hyperfinite, but there are hyperfinite sets
which are not finite. For example the set%
\begin{equation*}
\mathbb{R}^{\circ }:=\ \underset{\lambda \uparrow {\mathbb{U}}}{\lim }(%
\mathbb{R}\cap \lambda )
\end{equation*}%
is not finite. The hyperfinite sets are very important since they inherit
many properties of finite sets via Leibnitz principle. For example, $\mathbb{%
R}^{\circ }$ has the maximum and the minimum and every internal function%
\begin{equation*}
f:\mathbb{R}^{\circ }\rightarrow \mathbb{R}^{\ast }
\end{equation*}%
has the maximum and the minimum as well.

Also, it is possible to add the elements of an hyperfinite set of numbers or
vectors as follows: let%
\begin{equation*}
A:=\ \underset{\lambda \uparrow {\mathbb{U}}}{\lim }A_{\lambda }
\end{equation*}%
be an hyperfinite set; then the hyperfinite sum is defined in the following
way: 
\begin{equation*}
\sum_{a\in A}a=\ \underset{\lambda \uparrow {\mathbb{U}}}{\lim }\sum_{a\in
A_{\lambda }}a.
\end{equation*}%
In particular, if $A_{\lambda }=\left\{ a_{1}(\lambda ),...,a_{\beta
(\lambda )}(\lambda )\right\} \ $with\ $\beta (\lambda )\in \mathbb{N},\ $%
then setting 
\begin{equation*}
\beta =\ \underset{\lambda \uparrow {\mathbb{U}}}{\lim }\ \beta (\lambda
)\in \mathbb{N}^{\ast }
\end{equation*}%
we use the notation%
\begin{equation*}
\sum_{j=1}^{\beta }a_{j}=\ \underset{\lambda \uparrow {\mathbb{U}}}{\lim }%
\sum_{j=1}^{\beta (\lambda )}a_{j}(\lambda ).
\end{equation*}

\subsection{Qualified sets}

When we have a net $\varphi :Q\rightarrow ${$\mathbb{U}$}$_{n}$, where $%
Q\subset \mathcal{F}$, we can define the $\Lambda $-limit of $\varphi $ by
posing%
\begin{equation*}
\lim_{\lambda \in Q}\varphi (\lambda )=\lim_{\lambda \uparrow {\mathbb{U}}}%
\widetilde{\varphi }(\lambda )
\end{equation*}%
where 
\begin{equation*}
\widetilde{\varphi }(\lambda )=\left\{ 
\begin{array}{cc}
\varphi (\lambda ) & \text{for}\ \ \lambda \in Q \\ 
\varnothing & \text{for\ }\ \lambda \notin Q%
\end{array}%
\right.
\end{equation*}%
As one can expect, if two nets $\varphi ,\psi $ are equal on a "large" or a
"qualified" subset of $\mathcal{F}$ then they share the same $\Lambda $%
-limit. The notion of "qualified" subset of $\mathcal{F}$ can be precisely
defined as follows:

\begin{definition}
\label{qua}We say that a set $Q\subset \mathcal{F}$ is qualified if for
every bounded net $\varphi $ we have that 
\begin{equation*}
\lim_{\lambda \uparrow {\mathbb{U}}}\varphi (\lambda )=\lim_{\lambda \in
Q}\varphi (\lambda ).
\end{equation*}
\end{definition}

By the above definition, we have that the $\Lambda $-limit of a net $\varphi 
$ depends only on the values that $\varphi $ takes on a qualified set (it is
in this sense that we could imagine $Q$ to be "large"). It is easy to see
that (nontrivial) qualified sets exist. For example by ($\Lambda $-2) we
deduce that, for every $\lambda _{0}\in \mathcal{F}$, the set%
\begin{equation*}
Q\left( \lambda _{0}\right) :=\left\{ \lambda \in \mathcal{F}\ |\ \lambda
_{0}\subseteq \lambda \right\}
\end{equation*}%
is qualified. In this paper, we will use the notion of qualified set via the
following Theorem:

\begin{theorem}
\label{billo}Let $\mathcal{R}$ be a relation in {$\mathbb{U}$}$_{n}$ for
some $n\geq 0$ and let $\varphi $, $\psi :\mathcal{F\rightarrow \mathbb{U}}%
_{n}$. Then the following statements are equivalent:

\begin{itemize}
\item there exists a qualified set $Q$ such that 
\begin{equation*}
\forall \lambda \in Q,\ \varphi (\lambda )\mathcal{R}\psi (\lambda );
\end{equation*}

\item we have 
\begin{equation*}
\left( \underset{\lambda \uparrow {\mathbb{U}}}{\lim }\varphi (\lambda
)\right) \mathcal{R}^{\ast }\left( \underset{\lambda \uparrow {\mathbb{U}}}{%
\lim }\psi (\lambda )\right) .
\end{equation*}
\end{itemize}
\end{theorem}

\textbf{Proof}: It is an immediate consequence of Theorem \ref{limit} and
the definition of qualified set.

$\square $

\section{Ultrafunctions}

In this section, we will introduce the notion of ultrafunction and we will
analyze its first properties.

\subsection{Definition of Ultrafunctions}

Let $\Omega $ be a set in $\mathbb{R}^{N}$, and let $V(\Omega )\ $be a (real
or complex) vector space such that $\mathcal{D}(\overline{\Omega })\subseteq
V(\Omega )\subseteq L^{2}(\Omega )\cap \mathcal{C}(\overline{\Omega }).$

\begin{definition}
Given the function space $V(\Omega )$ we set%
\begin{equation*}
\widetilde{V}(\Omega ):=\lim_{\lambda \uparrow {\mathbb{U}}}V_{\lambda
}(\Omega )=Span^{\ast }(V(\Omega )^{\circ }),
\end{equation*}

where

\begin{equation*}
V_{\lambda }(\Omega )=Span(V(\Omega )\cap \lambda ).
\end{equation*}%
$\widetilde{V}(\Omega )$ will be called the \textbf{space of ultrafunctions}
generated by $V(\Omega ).$
\end{definition}

So, given any vector space of functions $V(\Omega )$, the space of
ultrafunction generated by $V(\Omega )$ is a vector space of hyperfinite
dimension that includes $V(\Omega )$, and the ultrafunctions are $\Lambda $%
-limits of functions in $V_{\lambda }$. Hence the ultrafunctions are
particular internal functions 
\begin{equation*}
u:\left( \mathbb{R}^{\ast }\right) ^{N}\rightarrow {\mathbb{C}^{\ast }.}
\end{equation*}

Observe that, by definition, the dimension of $\widetilde{V}(\Omega )$ (that
we denote by $\beta )$ is equal to the internal cardinality of any of its
bases, and the following formula holds: 
\begin{equation*}
\beta =\lim_{\lambda \uparrow {\mathbb{U}}}\text{dim}(V_{\lambda }(\Omega )).
\end{equation*}

Since $\widetilde{V}(\Omega )\subset \left[ L^{2}(\mathbb{R})\right] ^{\ast
},$ it can be equipped with the following scalar product%
\begin{equation*}
\left( u,v\right) =\int^{\ast }u(x)\overline{v(x)}\ dx,
\end{equation*}%
where $\int^{\ast }$ is the natural extension of the Lebesgue integral
considered as a functional%
\begin{equation*}
\int :L^{1}(\Omega )\rightarrow {\mathbb{C}}.
\end{equation*}
Notice that the Euclidean structure of $\widetilde{V}(\Omega )$ is the $%
\Lambda $-limit of the Euclidean structure of every $V_{\lambda }$ given by
the usual $L^{2}$ scalar product. The norm of an ultrafunction will be given
by 
\begin{equation*}
\left\Vert u\right\Vert =\left( \int^{\ast }|u(x)|^{2}\ dx\right) ^{\frac{1}{%
2}}.
\end{equation*}

\begin{remark}
\label{nina}Notice that the natural extension $f^{\ast }$ of a function $f$
is an ultrafunction if and only if $f\in V(\Omega ).$
\end{remark}

\textbf{Proof. }Let $f\in V(\Omega ),$ and $Q(f)=\{\lambda \in \mathcal{F}$ $%
|$ $f\in \lambda \}.$ Since, for every $\lambda $ $\mathbb{\in }$ $Q(f),$ $%
f\in V_{\lambda }(\Omega )$ and, as we observed in section \ref{qualified}, $%
Q(f)$ is a qualified set, it follows by Theorem \ref{billo} that $f^{\ast
}\in \widetilde{V}(\Omega ).$

Conversely, if $f\notin V(\Omega )$ then by Leibnitz Principle it follows
that $f^{\ast }\notin V^{\ast }(\Omega )$ and, since $\widetilde{V}(\Omega
)\subset V^{\ast }(\Omega )$, this entails the thesis.

$\square $

\subsection{Delta, Sigma and Theta Basis}

In this section we introduce three particular kinds of bases for $V\left(
\Omega \right) $ and we study their main properties. We start by defining
the \textit{Delta ultrafunctions}:

\begin{definition}
\label{dede}Given a number $q\in \Omega ^{\ast },$ we denote by $\delta
_{q}(x)$ an ultrafunction in $\widetilde{V}(\Omega )$ such that 
\begin{equation}
\forall v\in \widetilde{V}(\Omega ),\ \int^{\ast }v(x)\delta _{q}(x)dx=v(q).
\label{deltafunction}
\end{equation}%
$\delta _{q}(x)$ is called Delta (or the Dirac) ultrafunction concentrated
in $q$.
\end{definition}

Let us see the main properties of the Delta ultrafunctions:

\begin{theorem}
\label{delta} We have the following properties:

\begin{enumerate}
\item For every $q\in \Omega ^{\ast }$ there exists an unique Delta
ultrafunction concentrated in $q;$

\item for every $a,\ b\in \Omega ^{\ast }\ \delta _{a}(b)=\delta _{b}(a);$

\item $\left\Vert \delta _{q}\right\Vert ^{2}=\delta _{q}(q).$
\end{enumerate}
\end{theorem}

\textbf{Proof.} 1. Let $\left\{ e_{j}\right\} _{j=1}^{\beta }$ be an
orthonormal real basis of $\widetilde{V}(\Omega ),$ and set 
\begin{equation*}
\delta _{q}(x)=\sum_{j=1}^{\beta }e_{j}(q)e_{j}(x).\text{\label{deltaserie}}
\end{equation*}

Let us prove that $\delta _{q}(x)$ actually satisfies (\ref{deltafunction}).
Let $v(x)=\sum_{j=1}^{\beta }v_{j}e_{j}(x)$ be any ultrafunction. Then%
\begin{eqnarray*}
\int^{\ast }v(x)\delta _{q}(x)dx &=&\int^{\ast }\left( \sum_{j=1}^{\beta
}v_{j}e_{j}(x)\right) \left( \sum_{k=1}^{\beta }e_{k}(q)e_{k}(x)\right) dx=
\\
&=&\sum_{j=1}^{\beta }\sum_{k=1}^{\beta }v_{j}e_{k}(q)\int^{\ast
}e_{j}(x)e_{k}(x)dx= \\
&=&\sum_{j=1}^{\beta }\sum_{k=1}^{\beta }v_{j}e_{k}(q)\delta
_{j,q}=\sum_{j=1}^{\beta }v_{k}e_{k}(q)=v(q).
\end{eqnarray*}

So $\delta _{q}(x)$ is a Delta ultrafunction centered in $q$.

It is unique: if $f_{q}(x)$ is another Delta ultrafunction centered in $q$
then for every $y\in \Omega ^{\ast }$ we have:%
\begin{equation*}
\delta _{q}(y)-f_{q}(y)=\int^{\ast }(\delta _{q}(x)-f_{q}(x))\delta
_{y}(x)dx=\delta _{y}(q)-\delta _{y}(q)=0
\end{equation*}

and hence $\delta _{q}(y)=f_{q}(y)$ for every $y\in \Omega ^{\ast }.$

2.$\ \delta _{a}\left( b\right) =\int^{\ast }\delta _{a}(x)\delta _{b}(x)\
dx=\delta _{b}\left( a\right) .$

3. $\left\Vert \delta _{q}\right\Vert ^{2}=\int^{\ast }\delta _{q}(x)\delta
_{q}(x)=\delta _{q}(q)$.\newline

$\square $

\begin{definition}
A Delta-basis $\left\{ \delta _{a}(x)\right\} _{a\in \Sigma }$ $(\Sigma
\subset \Omega ^{\ast })$ is a basis for $\widetilde{V}(\Omega )$ whose
elements are Delta ultrafunctions. Its dual basis $\left\{ \sigma
_{a}(x)\right\} _{a\in \Sigma }$ is called Sigma-basis. We recall that, by
definition of dual basis, for every $a,b\in \Omega ^{\ast }$ the equation%
\begin{equation}
\int^{\ast }\delta _{a}(x)\sigma _{b}(x)dx=\delta _{ab}  \label{mimma}
\end{equation}%
holds. The set $\Sigma \subset \Omega ^{\ast }$ is called set of independent
points.
\end{definition}

The existence of a Delta-basis is an immediate consequence of the following
fact:

\begin{remark}
The set $\left\{ \delta _{a}(x)|a\in \Omega ^{\ast }\right\} $ generates all 
$\widetilde{V}(\Omega ).$ In fact, let $G(\Omega )$ be the vectorial space
generated by the set $\left\{ \delta _{a}(x)\ |\ a\in \Omega ^{\ast
}\right\} $ and suppose that $G(\Omega )$ is properly included in $%
\widetilde{V}(\Omega ).$ Then the orthogonal $G(\Omega )^{\perp }$ of $%
G(\Omega )$ in $\widetilde{V}(\Omega )$ contains a function $f\neq 0.$ But,
since $f\in $ $G(\Omega )^{\perp },$ for every $a\in \Omega ^{\ast \text{ }}$%
we have 
\begin{equation*}
f(a)=\int^{\ast }f(x)\delta _{a}(x)dx=0,
\end{equation*}%
so $f_{\upharpoonleft _{\Omega ^{\ast }}}=0$ and this is absurd. Thus the
set $\left\{ \delta _{a}(x)\ |\ a\in \Omega ^{\ast }\right\} $ generates $%
\widetilde{V}(\Omega ),$ hence it contains a basis.
\end{remark}

Let us see some properties of Delta- and Sigma-bases:

\begin{theorem}
\label{tbase}A Delta-basis $\left\{ \delta _{q}(x)\right\} _{q\in \Sigma }$
and its dual basis $\left\{ \sigma _{q}(x)\right\} _{q\in \Sigma }$ satisfy
the following properties:

\begin{enumerate}
\item if $u\in \widetilde{V}(\Omega )$, then%
\begin{equation*}
u(x)=\sum_{q\in \Sigma }\left( \int^{\ast }\sigma _{q}(\xi )u(\xi )d\xi
\right) \delta _{q}(x);
\end{equation*}

\item if $u\in \widetilde{V}(\Omega )$, then%
\begin{equation}
u(x)=\sum_{q\in \Sigma }u(q)\sigma _{q}(x);  \label{brava+}
\end{equation}

\item if two ultrafunctions $u$ and $v$ coincide on a set of independent
points then they are equal;

\item if $\Sigma $ is a set of independent points and $a,b\in \Sigma $ then $%
\sigma _{a}(b)=\delta _{ab};$

\item for any $q\in \Omega ^{\ast },$ $\sigma _{q}(x)$ is well defined.
\end{enumerate}
\end{theorem}

\textbf{Proof. }1. It is an immediate consequence of the definition of dual
basis.

2. Since $\left\{ \delta _{q}(x)\right\} _{q\in \Sigma }$ is the dual basis
of $\left\{ \sigma _{q}(x)\right\} _{q\in \Sigma }$ we have that%
\begin{equation*}
u(x)=\sum_{q\in \Sigma }\left( \int \delta _{q}(\xi )u(\xi )d\xi \right)
\sigma _{q}(x)=\sum_{q\in \Sigma }u(q)\sigma _{q}(x).
\end{equation*}

3. It follows directly from 2.

4. If follows directly by equation (\ref{mimma})

5. Given any point $q\in \Omega ^{\ast }$ clearly there is a Delta-basis $%
\left\{ \delta _{a}(x)\right\} _{a\in \Sigma }$ with $q\in \Sigma .$ Then $%
\sigma _{q}(x)$ can be defined by mean of the basis $\left\{ \delta
_{a}(x)\right\} _{a\in \Sigma }.$ We have to prove that, given another Delta
basis $\left\{ \delta _{a}(x)\right\} _{a\in \Sigma ^{\prime }}$ with $q\in
\Sigma ^{\prime },$ the corresponding $\sigma _{q}^{\prime }(x)$ is equal to 
$\sigma _{q}(x).$ Using (2), with $u(x)=\sigma _{q}^{\prime }(x),$ we have
that%
\begin{equation*}
\sigma _{q}^{\prime }(x)=\sum_{a\in \Sigma }\sigma _{q}^{\prime }(a)\sigma
_{a}(x).
\end{equation*}%
Then, by (4), it follows that $\sigma _{q}^{\prime }(x)=\sigma _{q}(x).$

$\square $

Let $\Sigma \mathbf{\ }$be a set of independent points, and let $L_{\Sigma }:%
\widetilde{V}(\Omega )\rightarrow \widetilde{V}(\Omega )$ be the linear
operator such that%
\begin{equation*}
L_{\Sigma }\sigma _{a}(x)=\delta _{a}(x)
\end{equation*}%
for every $a\in \Sigma .$

\begin{proposition}
$L_{\Sigma }$ is selfadjoint, positive and%
\begin{equation*}
\int^{\ast }L_{\Sigma }u(x)v(x)dx=\sum_{a\in \Sigma }u(a)v(a).
\end{equation*}
\end{proposition}

\textbf{Proof. }Since $u(x)=\sum_{a\in \Sigma }u(a)\sigma _{a}(x)$ and $%
v(x)=\sum_{a\in \Sigma }v(a)\sigma _{a}(x),$ then%
\begin{eqnarray*}
\int^{\ast }L_{\Sigma }u(x)v(x)dx &=&\int^{\ast }L_{\Sigma }\left(
\sum_{a\in \Sigma }u(a)\sigma _{a}(x)\right) \left( \sum_{b\in \Sigma
}v(b)\sigma _{b}(x)\right) dx= \\
&=&\sum_{a\in \Sigma }\sum_{b\in \Sigma }u(a)v(b)\int^{\ast }\delta
_{a}(x)\sigma _{b}(x)dx=\sum_{a\in \Sigma }u(a)v(a).
\end{eqnarray*}%
Hence, clearly, $L_{\Sigma }$\textbf{\ }is selfadjoint and positive.

$\square $

From now on, we consider the set $\Sigma \mathbf{\ }$fixed once for all and
we simply denote the operator $L_{\Sigma }$ by $L.$ Since $L$ is a positive
selfadjoint operator, $A=L^{1/2}$ is a well defined positive selfadjoint
operator. For every $a\in \Sigma $ we set 
\begin{equation*}
\theta _{a}(x)=A\sigma _{a}(x).
\end{equation*}

\begin{theorem}
\label{burro}The following properties hold:

\begin{enumerate}
\item $\left\{ \theta _{a}(x)\right\} _{a\in \Sigma }$ is an orthonormal
basis;

\item for every $a,b\in \Sigma ,$ $\theta _{a}(b)=\theta _{b}(a)$;

\item for every ultrafunction $u$ we have%
\begin{equation*}
u(x)=\sum_{a\in \Sigma }u(a)\sigma _{a}(x)=\sum_{a\in \Sigma }\underline{u}%
(a)\theta _{a}(x)=\sum_{a\in \Sigma }\underline{\underline{u}}(a)\delta
_{a}(x),
\end{equation*}%
where we have set, for every $a\in \Sigma ,$%
\begin{equation*}
\underline{u}(a):=(A^{-1}u)(a)=\int^{\ast }\theta _{a}(\xi )u(\xi )d\xi ;
\end{equation*}%
\begin{equation*}
\ \ \underline{\underline{u}}(a)=(A^{-1}\underline{u})(a)=(L^{-1}u)(a)=%
\int^{\ast }\sigma _{a}(\xi )u(\xi )d\xi ;
\end{equation*}

\item for every ultrafunctions $u,v$ we have%
\begin{equation*}
\int^{\ast }u(x)v(x)dx=\sum_{a\in \Sigma }\underline{u}(a)\underline{v}%
(a)=\sum_{a\in \Sigma }\underline{\underline{u}}(a)v(a);
\end{equation*}

\item for every ultrafunction $u$\ we have%
\begin{equation*}
\int^{\ast }u(x)dx=\sum_{a\in \Sigma }\underline{\underline{u}}(a).
\end{equation*}
\end{enumerate}
\end{theorem}

\textbf{Proof: }1) $\left\{ \theta _{a}(x)\right\} _{a\in \Sigma }$ is a
basis since it is the image of the basis $\left\{ \sigma _{a}(x)\right\}
_{a\in \Sigma }$ respect to the invertible linear application $L.$ It is
orthonormal: for every $a,b\in \Sigma $ we have%
\begin{eqnarray*}
\int^{\ast }\theta _{a}(x)\theta _{b}(x)dx &=&\int^{\ast }A\sigma
_{a}(x)A\sigma _{b}(x)dx=\int^{\ast }L\sigma _{a}(x)\sigma _{b}(x)= \\
&=&\sigma _{b}(a)=\delta _{ab.}
\end{eqnarray*}

2) We have%
\begin{eqnarray*}
\theta _{a}(b) &=&\int^{\ast }\theta _{a}(x)\delta _{b}(x)dx=\int^{\ast
}\theta _{a}(x)A\theta _{b}(x)dx= \\
&=&\int^{\ast }A\theta _{a}(x)\theta _{b}(x)dx=\int^{\ast }\delta
_{a}(x)\theta _{b}(x)dx=\theta _{b}(a).
\end{eqnarray*}

3) The equality%
\begin{equation*}
u(x)=\sum_{a\in \Sigma }u(a)\sigma _{a}(x)
\end{equation*}%
has been proved in Theorem \ref{tbase}, (\ref{brava+}); the equality 
\begin{equation*}
u(x)=\sum_{a\in \Sigma }\underline{u}(a)\theta _{a}(x),
\end{equation*}%
where $\underline{u}(a)=$ $\int^{\ast }\theta _{a}(\xi )u(\xi )d\xi $,
follows since $\left\{ \theta _{a}(x)\right\} _{a\in \Sigma }$ is an
orthonormal basis. And 
\begin{eqnarray*}
(A^{-1}u)(a) &=&\int^{\ast }\delta _{a}(\xi )A^{-1}u(\xi )d\xi = \\
&=&\int^{\ast }A^{-1}\delta _{a}(\xi )u(\xi )d\xi =\int^{\ast }\theta
_{a}(\xi )u(\xi )d\xi
\end{eqnarray*}%
since $A$ (and, so, $A^{-1})$ is selfadjoint.

The equality 
\begin{equation*}
u(x)=\sum_{a\in \Sigma }\underline{\underline{u}}(a)\delta _{a}(x),
\end{equation*}%
where $\ \underline{\underline{u}}(a)=\int^{\ast }\sigma _{a}(\xi )u(\xi
)d\xi $, follows by point (1) in Theorem \ref{tbase}. And $\ $%
\begin{eqnarray*}
\underline{\underline{u}}(a) &=&\int^{\ast }\sigma _{a}(\xi )u(\xi )d\xi
=\int^{\ast }L^{-1}\delta _{a}(\xi )u(\xi )d\xi = \\
&=&\int^{\ast }\delta _{a}(\xi )L^{-1}u(\xi )d\xi =(L^{-1}u)(a).
\end{eqnarray*}

4) We have that $\int^{\ast }u(x)v(x)dx=\sum_{a\in \Sigma }\underline{u}(a)%
\underline{v}(a)$ since $\left\{ \theta _{a}(x)\right\} _{a\in \Sigma }$ is
an orthonormal basis:%
\begin{eqnarray*}
\int^{\ast }u(x)v(x)dx &=&\int^{\ast }\left( \sum_{a\in \Sigma }\underline{u}%
(a)\theta _{a}(x)\right) \left( \sum_{b\in \Sigma }\underline{v}(b)\theta
_{b}(x)dx\right) = \\
&=&\sum_{a\in \Sigma }\sum_{b\in \Sigma }\underline{u}(a)\underline{v}%
(b)\int^{\ast }\theta _{a}(x)\theta _{b}(x)dx=\sum_{a\in \Sigma }\underline{u%
}(a)\underline{v}(a);
\end{eqnarray*}

the equality $\int^{\ast }u(x)v(x)dx=\sum_{a\in \Sigma }\underline{%
\underline{u}}(a)v(a)$ follows by expressing $u(x)$ in the Delta basis and $%
v(x)$ in the Sigma basis:%
\begin{eqnarray*}
\int^{\ast }u(x)v(x)dx &=&\int^{\ast }\left( \sum_{a\in \Sigma }\underline{%
\underline{u}}(a)\delta _{a}(x)\right) \left( \sum_{b\in \Sigma }v(b)\sigma
_{b}(x)\right) dx= \\
&=&\sum_{a\in \Sigma }\sum_{b\in \Sigma }v(b)\underline{\underline{u}}%
(a)\int^{\ast }\delta _{a}(x)\sigma _{b}(x)dx=\sum_{a\in \Sigma }\underline{%
\underline{u}}(a)v(a).
\end{eqnarray*}

5) This follows by expressing $u(x)$ in the Delta basis: 
\begin{equation*}
\int^{\ast }u(x)dx=\int^{\ast }\sum_{a\in \Sigma }\underline{\underline{u}}%
(a)\delta _{a}(x)dx=\sum_{a\in \Sigma }\underline{\underline{u}}%
(a)\int^{\ast }\delta _{a}(x)dx=\sum_{a\in \Sigma }\underline{\underline{u}}%
(a).
\end{equation*}

$\square $

\subsection{Canonical extension of a function}

Let $V^{\prime }(\Omega )$ denote the dual of $V(\Omega )\ $and let $%
\mathfrak{M}$ denote the set of measurable functions in $\mathbb{R}^{N}$. If 
$T\in V^{\prime }(\Omega )$ and if there is a function $f\in \mathfrak{M}$
such that 
\begin{equation*}
\forall v\in V(\Omega ),\ \left\langle T,v\right\rangle =\int f(x)v(x)dx
\end{equation*}%
then $T\ $and $f$ will be identified, and with some abuse of notation we
shall write $T=f\in V^{\prime }(\Omega )\cap \mathfrak{M.}$ With this
identification, $V^{\prime }(\Omega )\cap \mathfrak{M}\subset L^{2}.$

\begin{definition}
\label{tilde} If $T\in \left[ V^{\prime }(\Omega )\right] ^{\ast },$ there
exists a unique ultrafunction$\widetilde{\ T}(x)$ such that 
\begin{equation*}
\forall v\in \widetilde{V}(\Omega ),\ \left\langle T,v\right\rangle
=\int^{\ast }\widetilde{T}(x)v(x)dx.
\end{equation*}%
In particular, if $u\in \left[ V^{\prime }(\Omega )\cap \mathfrak{M}\right]
^{\ast }\mathfrak{,}$ $\widetilde{u}$ will denote the unique ultrafunction
such that%
\begin{equation*}
\forall v\in \widetilde{V}(\Omega ),\ \int^{\ast }u(x)v(x)dx=\int^{\ast }%
\widetilde{u}(x)v(x)dx.
\end{equation*}
\end{definition}

Notice that $V^{\prime }(\Omega )\cap \mathfrak{M}$ is a space of
distributions which contains the delta measures, so to every Delta
distribution $\delta _{q}$ is associated an ultrafunction which, by
definition, is the Delta ultrafunction centered in $q$, as expected.

\begin{definition}
If $f\in V^{\prime }(\Omega )\cap \mathfrak{M,}$ $\widetilde{(f^{\ast })}$
is called the {canonical extension of $f$. In the following, since }$f$ and{%
\ }$f^{\ast }$ can be identified, we will write $\widetilde{f}$ instead of $%
\widetilde{(f^{\ast })}.$
\end{definition}

Thus any function%
\begin{equation*}
f:\mathbb{R}^{N}\rightarrow \mathbb{R}
\end{equation*}%
can be extended to the function%
\begin{equation*}
f\mathbb{^{\ast }}:(\mathbb{R^{\ast })}^{N}\rightarrow \mathbb{R^{\ast }}
\end{equation*}%
which is called the natural extension of $f$ and if $f\in V^{\prime }(\Omega
)\cap \mathfrak{M,}$ we have also the canonical extension of $f$ given by 
\begin{equation*}
\widetilde{f}:(\mathbb{R^{\ast })}^{N}\rightarrow \mathbb{R^{\ast }}
\end{equation*}%
If $f\notin V(\Omega ),\ $by Remark \ref{nina}, $\widetilde{f}\neq f^{\ast }$%
, thus$\ f^{\ast }\notin \widetilde{V}(\Omega ).$

\bigskip

\textbf{Example:} if $\Omega =(-1,1),\ $then $|x|^{-1/2}\in V(-1,1)^{\prime
}\cap \mathfrak{M};$ the ultrafunction $\widetilde{|x|^{-1/2}}$ is different
from $\left( |x|^{-1/2}\right) ^{\ast }$ since the latter is not defined for 
$x=0,$ while 
\begin{equation*}
\left( \widetilde{|x|^{-1/2}}\right) _{x=0}=\int^{\ast }|x|^{-1/2}\delta
_{0}(x)dx.
\end{equation*}

\begin{theorem}
\label{CA} If $T\in \left[ V(\Omega )^{\prime }\right] ^{\ast },$ then 
\begin{eqnarray*}
\widetilde{T}(x) &=&\sum_{q\in \Sigma }\left\langle T,\delta
_{q}\right\rangle \sigma _{q}(x)= \\
&=&\sum_{q\in \Sigma }\left\langle T,\theta _{q}\right\rangle \theta _{q}(x)=
\\
&=&\sum_{q\in \Sigma }\left\langle T,\sigma _{q}\right\rangle \delta _{q}(x).
\end{eqnarray*}%
In particular, if $f\in \left[ V^{\prime }(\Omega )\cap \mathfrak{M}\right]
^{\ast }$ 
\begin{eqnarray}
\widetilde{f}(x) &=&\sum_{q\in \Sigma }\left[ \int f^{\ast }(\xi )\delta
_{q}(\xi )d\xi \right] \sigma _{q}(x)=  \label{bella} \\
&=&\sum_{q\in \Sigma }\left[ \int f^{\ast }(\xi )\theta _{q}(\xi )d\xi %
\right] \theta _{q}(x)= \\
&=&\sum_{q\in \Sigma }\left[ \int f^{\ast }(\xi )\sigma _{q}(\xi )d\xi %
\right] \delta _{q}(x).
\end{eqnarray}
\end{theorem}

\textbf{Proof.} It is sufficient to prove that 
\begin{equation*}
\forall v\in V(\Omega ),\ \int \sum_{q\in \Sigma }\left\langle T,\delta
_{q}\right\rangle \sigma _{q}(x)v(x)dx=\left\langle T,v\right\rangle .
\end{equation*}%
We have that 
\begin{eqnarray*}
\ \int \sum_{q\in \Sigma }\left\langle T,\delta _{q}\right\rangle \sigma
_{q}(x)v(x)dx &=&\sum_{q\in \Sigma }\left\langle T,\delta _{q}\right\rangle
\int \sigma _{q}(x)v(x)dx= \\
&=&\left\langle T,\sum_{q\in \Sigma }\left( \int \sigma _{q}(x)v(x)dx\right)
\delta _{q}\right\rangle =\left\langle T,v\right\rangle .
\end{eqnarray*}

The other equalities can be proved similarly.

$\square $

\subsection{Ultrafunctions and distributions}

In this section we will show that the space of ultrafunctions is reacher
than the space of distribution, in the sense that any distribution can be
represented by an ultrafunction and that the converse is not true.

\begin{definition}
Let $D\subset \widetilde{V}(\Omega )$ be a vector space. We say that two
ultrafunctions $u$ and $v$ are $D$-equivalent if%
\begin{equation*}
\forall \varphi \in D,\ \int^{\ast }\left( u(x)-v(x)\right) \varphi (x)dx=0.
\end{equation*}%
We say that two ultrafunctions $u$ and $v$ are distributionally equivalent
if they are $\mathcal{D}(\Omega )$-equivalent.
\end{definition}

\begin{theorem}
Given $T\in \mathcal{D}^{\prime },$ there exists an ultrafunction $u$ such
that%
\begin{equation}
\forall \varphi \in \mathcal{D}(\Omega ),\ \int^{\ast }u(x)\varphi ^{\ast
}(x)dx=\left\langle T,\varphi \right\rangle .  \label{pippo}
\end{equation}
\end{theorem}

\textbf{Proof:} Let $\left\{ e_{j}(x)\right\} _{j\in J}$ be an orthonormal
basis of the hyperfinite space $\widetilde{V}(\Omega )\cap \mathcal{D}%
(\Omega )^{\ast }$ and take 
\begin{equation*}
u(x)=\sum_{j\in J}\left\langle T^{\ast },e_{j}\right\rangle \ e_{j}(x).
\end{equation*}

Now take $\varphi \in \mathcal{D}.$ Since $\varphi ^{\ast }\in \widetilde{V}%
(\Omega )\cap \mathcal{D}(\Omega )^{\ast },$ we have that%
\begin{equation*}
\varphi ^{\ast }(x)=\sum_{j\in J}\left( \int^{\ast }\varphi ^{\ast }(\xi
)e_{j}(\xi )d\xi \right) e_{j}(x).
\end{equation*}%
Thus 
\begin{eqnarray*}
\int^{\ast }u(x)\varphi ^{\ast }(x)dx &=&\int^{\ast }\sum_{j\in
J}\left\langle T^{\ast },e_{j}\right\rangle \ e_{j}(x)\varphi ^{\ast
}(x)dx=\sum_{j\in J}\left\langle T^{\ast },e_{j}\int^{\ast }e_{j}(x)\varphi
^{\ast }(x)dx\right\rangle = \\
&=&\left\langle T^{\ast },\sum_{j\in J}\left( \int^{\ast }e_{j}(x)\varphi
^{\ast }(x)dx\right) e_{j}\right\rangle =\left\langle T^{\ast },\varphi
^{\ast }\right\rangle =\left\langle T,\varphi \right\rangle .
\end{eqnarray*}

$\square $

\bigskip

The following proposition shows that the ultrafunction $u$ associated to the
distribution $T$ by (\ref{pippo}) is not unique:

\begin{proposition}
Take $T\in \mathcal{D}^{\prime }(\Omega )$ and let 
\begin{equation*}
V_{T}=\{u\in \widetilde{V}(\Omega ):\forall \varphi \in \mathcal{D}(\Omega
),\ \int^{\ast }u(x)\varphi ^{\ast }(x)dx=\left\langle T,\varphi
\right\rangle \},
\end{equation*}%
let $u\in V_{T}$ and let $v$ be any ultrafunction. Then

\begin{enumerate}
\item $v\in V_{T}$ if and only if $u$ and $v$ are $\mathcal{D}(\Omega )$%
-equivalent;

\item $V_{T}$ is infinite.
\end{enumerate}
\end{proposition}

\textbf{Proof:} 1) If $v\in V_{T}$ then $\forall \varphi \in \mathcal{D}%
(\Omega ),\ \int^{\ast }(u(x)-v(x))\varphi ^{\ast }(x)dx=\left\langle
T,\varphi \right\rangle -\left\langle T,\varphi \right\rangle =0$, so $u$
and $v$ are $\mathcal{D}(\Omega )$-equivalent; conversely, if $u$ and $v$
are $\mathcal{D}$-equivalent then $\forall \varphi \in \mathcal{D}(\Omega
),\ \int^{\ast }u(x)\varphi ^{\ast }(x)dx=\int^{\ast }v(x)\varphi ^{\ast
}(x)dx.$ Since $\int^{\ast }u(x)\varphi ^{\ast }(x)dx=\left\langle T,\varphi
\right\rangle $ then $v\in V_{T}.$

2) Let $v\neq 0$ be any ultrafunction in the orthogonal (in $\widetilde{V}%
(\Omega ))$ of $\widetilde{V}(\Omega )\cap \mathcal{D}(\Omega )^{\ast }$.
Then $u$ and $u+v$ are $\mathcal{D}(\Omega )$-equivalent, since $\int^{\ast
}(u(x)+v(x))\varphi ^{\ast }(x)dx=\int^{\ast }u(x)\varphi ^{\ast
}(x)dx+\int^{\ast }v(x)\varphi ^{\ast }(x)dx=\int^{\ast }u(x)\varphi ^{\ast
}(x)dx+0$. Since the orthogonal of $\widetilde{V}(\Omega )\cap \mathcal{D}%
(\Omega )^{\ast }$ is infinite, we obtain the thesis.

$\square $

\begin{remark}
There is a natural way to associate a unique ultrafunction to a distribution
(see also \cite{ultra}). In order to do this it is sufficient to split $%
\widetilde{V}(\Omega )$ in two orthogonal component: $\widetilde{V}(\Omega
)\cap \mathcal{D}(\Omega )^{\ast }$ and $\left( \widetilde{V}(\Omega )\cap 
\mathcal{D}(\Omega )^{\ast }\right) ^{\perp }.$ As we have seen in the proof
of the above theorem every ultrafunction in $V_{T}$ can be spitted in two
components, $u+v$ where $v\in \left( \widetilde{V}(\Omega )\cap \mathcal{D}%
(\Omega )^{\ast }\right) ^{\perp }\ $and $u\in \widetilde{V}(\Omega )\cap 
\mathcal{D}(\Omega )^{\ast }$ is univocally determined. Then, we have an
injective map%
\begin{equation*}
i:\mathcal{D}^{\prime }(\Omega )\rightarrow \widetilde{V}(\Omega )
\end{equation*}%
given by $i(T)=u$ where $u\in V_{T}\cap \mathcal{D}(\Omega )^{\ast }.$
\end{remark}

\begin{remark}
The space of ultrafunctions is richer than the space of distributions; for
example consider the function 
\begin{equation*}
u(x):=f(x)\min \left( x^{-2},\alpha \right) 
\end{equation*}%
where $\alpha >0$ is an infinite number and $f(x)$ is a function with
compact support such that $f(0)=1$. Since $u\in \left[ V^{\prime }(\Omega
)\cap \mathfrak{M}\right] ^{\ast }\mathfrak{,}$ $\widetilde{u}$ is well
defined (see def. \ref{tilde}). On the other hand, $\widetilde{u}$ does not
correspond to any distribution since 
\begin{equation*}
\int^{\ast }\widetilde{u}(x)\varphi ^{\ast }(x)dx=\int^{\ast }f^{\ast
}(x)\min \left( x^{-2},\alpha \right) \varphi ^{\ast }(x)dx
\end{equation*}%
is infinite when $\varphi (x)\geq 0$ and $\varphi (0)>0.$ In \cite{ultra}
Section 6, it is presented an elliptic problem which has a solution in the
space of ultrafunctions, but no solution in the space of distributions.
\end{remark}

\section{Operations with ultrafunctions}

\subsection{Extension of operators\label{EO}}

\begin{definition}
\label{CE}Given the operator $F:V\left( \Omega \right) \rightarrow \mathcal{D%
}^{\prime }\left( \Omega \right) ,$ the map%
\begin{equation*}
\widetilde{F}:\widetilde{V}\left( \Omega \right) \rightarrow \widetilde{V}%
\left( \Omega \right)
\end{equation*}%
defined by%
\begin{equation}
\widetilde{F}\left( u\right) =\widetilde{F^{\ast }\left( u\right) }
\end{equation}%
is called \textbf{canonical }extension of $F$ ("$\sim $" is defined by \ref%
{tilde}).
\end{definition}

By the definition of $\widetilde{F}$, we have that%
\begin{equation}
\forall v\in \widetilde{V}\left( \Omega \right) ,\ \int^{\ast }\widetilde{F}%
\left( u(x)\right) v(x)\ dx=\int^{\ast }F^{\ast }\left( u(x)\right) v(x)dx.
\label{bellina}
\end{equation}

Comparing Definition \ref{CE} with Theorem \ref{CA} we have that%
\begin{eqnarray*}
\widetilde{F}(u(x)) &=&\sum_{q\in \Sigma }\left\langle F^{\ast }\left(
u\right) ,\delta _{q}\right\rangle \sigma _{q}(x)= \\
&=&\sum_{q\in \Sigma }\left\langle F^{\ast }\left( u\right) ,\theta
_{q}\right\rangle \theta _{q}(x)= \\
&=&\sum_{q\in \Sigma }\left\langle F^{\ast }\left( u\right) ,\sigma
_{q}\right\rangle \delta _{q}(x).
\end{eqnarray*}%
In particular, if $F:V\left( \Omega \right) \rightarrow V^{\prime }(\Omega
)\cap \mathfrak{M}^{\ast }:$%
\begin{eqnarray}
\widetilde{F}(u(x)) &=&\sum_{q\in \Sigma }\left[ \int F^{\ast }\left( u(\xi
)\right) \delta _{q}(\xi )d\xi \right] \sigma _{q}(x)=  \label{lola} \\
&=&\sum_{q\in \Sigma }\left[ \int F^{\ast }\left( u(\xi )\right) \theta
_{q}(\xi )d\xi \right] \theta _{q}(x)=  \notag \\
&=&\sum_{q\in \Sigma }\left[ \int F^{\ast }\left( u(\xi )\right) \sigma
_{q}(\xi )d\xi \right] \delta _{q}(x).  \notag
\end{eqnarray}

\subsection{Derivative}

A good generating space to define the derivative of an ultrafunction is the
following one:%
\begin{equation*}
V^{1}(\Omega )=H^{1,1}(\Omega )\cap \mathcal{C}(\overline{\Omega })\subseteq
L^{2}(\Omega )\cap \mathcal{C}(\overline{\Omega }).
\end{equation*}

In order to simplify the exposition, we will assume that $\Omega \subseteq 
\mathbb{R}$. The generalization of the notions exposed in this section is
immediate.

Let $u\in \widetilde{V^{1}}(\Omega )$ be a ultrafunction. Since $%
V^{1}(\Omega )^{\ast }\subset H^{1}(\Omega )^{\ast },$ we have that the
derivative $\frac{du}{dx}=\partial u=u^{\prime }$ is in $L^{2}(\Omega
)\subset \left[ V_{G}^{\prime }\cap \mathfrak{M}\right] ^{\ast }.$ Then we
can apply Definition \ref{CE}:

\begin{definition}
We set 
\begin{equation*}
Du=\widetilde{\partial }u=\widetilde{\partial u}.
\end{equation*}%
The operator 
\begin{equation*}
D:\widetilde{V^{1}}(\Omega )\rightarrow \widetilde{V^{1}}(\Omega )
\end{equation*}%
is called (generalized) derivative of the ultrafunction $u.$
\end{definition}

By (\ref{lola}) we have the following representation of the derivative:%
\begin{equation*}
\forall u\in \widetilde{V^{1}}(\Omega ),\ Du(x)=\sum_{q\in \Sigma }\left[
\int^{\ast }u^{\prime }(\xi )\delta _{q}(\xi )d\xi \right] \sigma _{q}(x).
\end{equation*}

If $u^{\prime }\in \widetilde{V^{1}}(\Omega )\subset \left[ V^{1}(\Omega )%
\right] ^{\ast },$ we have that%
\begin{equation*}
Du(x)=\sum_{q\in \Sigma }u^{\prime }(q)\sigma _{q}(x)=u^{\prime }(x).
\end{equation*}%
In particular, if $u\in H^{2,1}(\Omega )\cap \mathcal{C}^{1}(\overline{%
\Omega }),$ $Du=u^{\prime }$ and so $D$ extends the operator $\frac{d}{dx}%
:H^{2,1}(\Omega )\cap \mathcal{C}^{1}(\overline{\Omega })\rightarrow
V^{1}(\Omega )$ to the operator $D:\widetilde{V^{1}}(\Omega )\rightarrow 
\widetilde{V^{1}}(\Omega ).$

\subsection{Fourier transform}

In this section we will investigate the extension of the one-dimensional
Fourier transform. A good space to work with the Fourier transform is the
space

\begin{equation*}
V^{\mathfrak{F}}(\mathbb{R})=H^{1}(\mathbb{R})\cap L^{2}(\mathbb{R},|x|^{2}).
\end{equation*}%
It is easy to see that the space $V^{\mathfrak{F}}(\mathbb{R})$ can be
characterized as follows:%
\begin{equation*}
V^{\mathfrak{F}}(\mathbb{R})=\left\{ u\in H^{1}(\mathbb{R}):\hat{u}\in H^{1}(%
\mathbb{R})\right\} .
\end{equation*}

In fact, if $\hat{u}\in H^{1}(\mathbb{R}),$ then $\int |\nabla u(\xi )|^{2}\
d\xi <+\infty $ and hence $\int |u(x)|^{2}|x|^{2}\ dx<+\infty .$

Then $V^{\mathfrak{F}}(\mathbb{R})\subset L^{2}(\mathbb{R},|x|^{2}),$ so $V^{%
\mathfrak{F}}(\mathbb{R})\subset H^{1}(\mathbb{R})\cap L^{2}(\mathbb{R}%
,|x|^{2})$ which is a Hilbert space equipped with the norm%
\begin{equation*}
\left\Vert u\right\Vert _{V^{\mathfrak{F}}(\mathbb{R})}^{2}=\int \left\vert
u(x)\right\vert ^{2}|x|^{2}dx+\int \left\vert \hat{u}(\xi )\right\vert
^{2}|\xi |^{2}d\xi .
\end{equation*}

Moreover 
\begin{eqnarray*}
\int \left\vert u(x)\right\vert dx &=&\int \left\vert u(x)\right\vert (1+|x|)%
\frac{1}{1+|x|}dx\leq \\
&\leq &\left( \int \left\vert u(x)\right\vert ^{2}(1+|x|)^{2}dx\right) ^{%
\frac{1}{2}}\left( \int \frac{1}{(1+|x|)^{2}}dx\right) ^{\frac{1}{2}}\leq \\
&\leq &const.\left( \left\Vert u\right\Vert _{L^{2}(\mathbb{R})}+\left\Vert
u\right\Vert _{L^{2}(\mathbb{R},|x|^{2})}\right) .
\end{eqnarray*}%
Thus, $V^{\mathfrak{F}}(\mathbb{R})\subset L^{1}(\mathbb{R}).$ Recalling
that the functions in $H^{1}(\mathbb{R})$ are continuous, we have that 
\begin{equation*}
V^{\mathfrak{F}}(\mathbb{R})\subset \mathcal{C}\left( \mathbb{R}\right) \cap
H^{1}(\mathbb{R})\cap L^{1}(\mathbb{R})\cap L^{2}(\mathbb{R},|x|^{2}).
\end{equation*}

We use the following definitions of Fourier transform: if $u\in \widetilde{%
V^{\mathfrak{F}}}(\mathbb{R})$, we set 
\begin{eqnarray}
\mathfrak{F}(u)(k) &=&\widehat{u}(k)=\frac{1}{\sqrt{2\pi }}\int^{\ast }u(x)\
e^{-ikx}\ dx;  \label{fur} \\
\mathfrak{F}^{-1}(u)(x) &=&\frac{1}{\sqrt{2\pi }}\int^{\ast }\widehat{u}(k)\
e^{ikx}\ dx.
\end{eqnarray}

Now, in order to deal with the Fourier transform in an easier way, we need a
new axiom whose consistency is easy to be verified (see Appendix):

\begin{axiom}
$(\mathbf{FTA)}$\textbf{(Fourier transform axiom)} If $u\in \widetilde{V^{%
\mathfrak{F}}}(\mathbb{R})$ then $\mathfrak{F}^{\ast }(u)\in \widetilde{V^{%
\mathfrak{F}}}(\mathbb{R})$ and $\bar{u}\in \widetilde{V^{\mathfrak{F}}}(%
\mathbb{R})$ (here $\bar{u}$ is the complex conjugate of $u).$
\end{axiom}

It is immediate to see that, by this axiom, for every ultrafunction, $u$ we
have 
\begin{equation*}
\mathfrak{F}^{\ast }(u)=\widetilde{\mathfrak{F}}(u)
\end{equation*}%
and hence, since there is no risk of ambiguity, we will simply write $%
\mathfrak{F}(u).$

It is well known that in the theory of tempered distributions we have that:%
\begin{equation*}
\mathfrak{F}(\delta _{a})=\frac{e^{-iak}}{\sqrt{2\pi }};
\end{equation*}%
\begin{equation*}
\mathfrak{F}\left( \frac{e^{iax}}{\sqrt{2\pi }}\right) =\delta _{a}.
\end{equation*}%
In the theory of ultrafunctions an analogous result holds:

\begin{proposition}
\label{marisa}We have that:

\begin{enumerate}
\item $\mathfrak{F}\left( \frac{\widetilde{e^{iax}}}{\sqrt{2\pi }}\right)
=\delta _{a}(k);$

\item $\mathfrak{F}\left( \delta _{a}(x)\right) =\frac{\widetilde{e^{-iak}}}{%
\sqrt{2\pi }};$

\item $\frac{1}{2\pi }\int^{\ast }\widetilde{e^{-iax}}\ \widetilde{e^{ikx}}\
dx=\delta _{a}(k).$
\end{enumerate}
\end{proposition}

\textbf{Proof.} 1. For every $v\in V^{\mathfrak{F}},$%
\begin{eqnarray*}
\int^{\ast }\mathfrak{F}\left( \frac{\widetilde{e^{iax}}}{\sqrt{2\pi }}%
\right) v(k)dk &=&\int^{\ast }\left( \frac{1}{2\pi }\int^{\ast }\widetilde{%
e^{-iak}}\ e^{ixk}dx\right) v(k)dk= \\
&=&\frac{1}{2\pi }\int^{\ast }\int^{\ast }\widetilde{e^{-iak}}\
e^{ixk}v(k)dkdx= \\
&=&\frac{1}{\sqrt{2\pi }}\int^{\ast }\widetilde{e^{-iak}}\mathfrak{F}%
^{-1}(v(k))dx=v(a).
\end{eqnarray*}

Hence, 1 holds.

2 - We have 
\begin{equation*}
\mathfrak{F}\left( \delta _{a}(x)\right) =\int^{\ast }\delta
_{a}(x)e^{-ikx}dx=\int^{\ast }\delta _{a}(x)\widetilde{e^{-ikx}}dx=%
\widetilde{e^{-ika}.}
\end{equation*}

3 - We have%
\begin{equation*}
\frac{1}{2\pi }\int^{\ast }\widetilde{e^{iax}}\ \widetilde{e^{-ikx}}\ dx=%
\frac{1}{2\pi }\int^{\ast }\widetilde{e^{iax}}\ e^{-ikx}dx=\mathfrak{F}%
\left( \frac{\widetilde{e^{iax}}}{\sqrt{2\pi }}\right) =\delta _{a}(k).
\end{equation*}

$\square $

\bigskip

By our definitions we have that:%
\begin{eqnarray*}
\widetilde{e^{ikx}} &=&\sum_{q\in \Sigma }\left[ \int^{\ast }e^{ik\xi
}\delta _{q}(\xi )d\xi \right] \sigma _{q}(x); \\
\widetilde{e^{ixk}} &=&\sum_{q\in \Sigma }\left[ \int^{\ast }e^{ix\xi
}\delta _{q}(\xi )d\xi \right] \sigma _{q}(k).
\end{eqnarray*}%
Therefore it is not evident whether $\widetilde{e^{ikx}}=\widetilde{e^{ixk}}$
or not. The following Corollary answers this question.

\begin{corollary}
We have that:%
\begin{equation*}
\widetilde{e^{ikx}}=\widetilde{e^{ixk}}.
\end{equation*}
\end{corollary}

\textbf{Proof.} By the previous proposition, we have that%
\begin{equation*}
\widetilde{e^{-ikx}}=\sqrt{2\pi }\mathfrak{F}\left( \delta _{k}(x)\right)
=\int^{\ast }\delta _{k}(x)e^{-ixk}dk=\int^{\ast }\delta _{x}(k)e^{-ixk}dx=%
\widetilde{e^{-ixk}}.
\end{equation*}%
Replacing $x$ with $-x$ we get the result.

$\square $

\bigskip

Since $\mathfrak{F:}V^{\mathfrak{F}}(\mathbb{R})\rightarrow V^{\mathfrak{F}}(%
\mathbb{R})$ is an isomorphism, it follows that, for any Delta-basis $%
\left\{ \delta _{a}\right\} _{a\in \Sigma },$ the set$\ $%
\begin{equation*}
\left\{ \frac{\widetilde{e^{iax}}}{\sqrt{2\pi }}\right\} _{a\in \Sigma
}=\left\{ \mathfrak{F}\left( \delta _{-a}\right) \right\} _{a\in \Sigma }
\end{equation*}%
is a basis and we get the following result:

\begin{theorem}
If $u\in V^{\mathfrak{F}}(\mathbb{R})$, then%
\begin{equation*}
u(x)=\frac{1}{\sqrt{2\pi }}\sum_{k\in \Sigma }\underline{\underline{\widehat{%
u}}}(k)\widetilde{e^{ikx}}.
\end{equation*}%
where we have set (see Theorem \ref{burro})%
\begin{equation*}
\underline{\underline{\widehat{u}}}(k)=\int^{\ast }\widehat{u}(\xi )\sigma
_{k}(\xi )d\xi .
\end{equation*}
\end{theorem}

\textbf{Proof. }Since $\left\{ \frac{\widetilde{e^{ikx}}}{\sqrt{2\pi }}%
\right\} _{k\in \Sigma }$ is a basis, any $u\in V^{\mathfrak{F}}(\mathbb{R})$
has the following representation:%
\begin{equation*}
u(x)=\frac{1}{\sqrt{2\pi }}\sum_{k\in \Sigma }u_{k}\widetilde{e^{ikx}.}
\end{equation*}

Let us compute the $u_{k}$'s: we have 
\begin{equation*}
\int \delta _{k}(x)\overline{\sigma _{b}(x)}dx=\int \delta _{k}(x)\sigma
_{b}(x)dx=\delta _{kb}
\end{equation*}%
and so 
\begin{equation*}
\int \widehat{\delta _{k}}(x)\overline{\widehat{\sigma _{b}}}(x)dx=\delta
_{kb}
\end{equation*}%
and by Proposition \ref{marisa}, 
\begin{equation*}
\int \frac{\widetilde{e^{-ikx}}}{\sqrt{2\pi }}\ \overline{\widehat{\sigma
_{b}}}(x)dx=\delta _{kb}.
\end{equation*}%
Hence $\left\{ \widehat{\sigma _{k}}(x)\right\} _{k\in \Sigma }\ $is the
dual basis of $\left\{ \frac{\widetilde{e^{-ikx}}}{\sqrt{2\pi }}\right\}
_{k\in \Sigma },$namely $\left\{ \widehat{\sigma _{k}}(-x)\right\} _{k\in
\Sigma }\ $is the dual basis of $\left\{ \frac{\widetilde{e^{ikx}}}{\sqrt{%
2\pi }}\right\} _{k\in \Sigma }.$ Hence, since $\widehat{\widehat{v}}%
(x)=v(-x),$ we have: 
\begin{eqnarray*}
u_{k} &=&\int u(\xi )\overline{\widehat{\sigma _{k}}(-\xi )}d\xi =\int u(\xi
)\overline{\widehat{\widehat{\sigma _{k}}}(-\xi )}d\xi = \\
&=&\int \widehat{u}(\xi )\overline{\sigma _{k}(\xi )}d\xi =\int \widehat{u}%
(\xi )\sigma _{k}(\xi )d\xi =\underline{\underline{\widehat{u}}}(k).
\end{eqnarray*}

$\square $

\section{Appendix\label{A}}

In this section we prove that the axiomatic construction of ultrafunctions
is coherent. We assume that the reader knows the key concepts in Nonstandard
Analysis (see e.g. \cite{keisler76}).

The following result has already been proved in \cite{ultra}. Here we give
an alternative proof of this result based on Nonstandard Analysis:

\begin{theorem}
\label{catullo}The set of axioms $\{(\Lambda $-$1)$\textsf{,}$(\Lambda $-$2)$%
,$(\Lambda $-$3)\}$ is consistent.
\end{theorem}

\textbf{Proof.}\ Let $\mathbb{U}$, $\mathbb{V}$ be mathematical universes
and let $\langle \mathbb{U},\mathbb{V},\star \rangle $ be a nonstandard
extension of $\mathbb{U}$ that is $|\mathbb{U}|^{+}$-saturated. We denote by 
$\mathcal{F}$ the set of finite subsets of $\mathbb{U}$ and, for every $%
\lambda \in $ $\mathcal{F}$, we pose%
\begin{equation*}
F_{\lambda }=\{S\subset \mathbb{V}|\text{ }S\text{ is hyperfinite and }%
\lambda ^{\star }\subset S\}.
\end{equation*}%
By saturation $\bigcap_{\lambda \in \mathcal{F}}F_{\lambda }\neq \emptyset .$
We take $\Lambda \in \bigcap_{\lambda \in \mathcal{F}}F_{\lambda }.$

For any given net $\varphi :\mathcal{F}\rightarrow \mathbb{U}$ we define its 
$\Lambda $-limit as%
\begin{equation*}
\lim_{\lambda \uparrow \mathbb{U}}\varphi (\lambda )=\varphi ^{\star
}(\Lambda )
\end{equation*}%
and we pose%
\begin{equation*}
\mathbb{K}=\lim\limits_{\lambda \uparrow \mathbb{U}}\mathbb{R}=\left\{
\lim\limits_{\lambda \uparrow \mathbb{U}}\varphi (\lambda )\mid \varphi :%
\mathcal{F}\rightarrow \mathbb{R}\right\} .
\end{equation*}

With these choices the $\Lambda $-limit satisfies the axioms ($\Lambda $-1)%
\textsf{,}($\Lambda $-2),($\Lambda $-3): the only nontrivial fact is ($%
\Lambda $-2). Let $\varphi $ be an eventually constant net, and let $\lambda
_{0}\in \mathcal{F},r\in \mathbb{R}$ be such that $\forall \lambda \in
\{\eta \in \mathcal{F}\mid \lambda _{0}\subset \eta \}$%
\begin{equation*}
\varphi (\lambda )=r.
\end{equation*}

By transfer it follows that$\ \forall \lambda \in \{\eta \in \mathcal{F}\mid
\lambda _{0}\subset \eta \}^{\star }=\{\eta \in \mathcal{F}^{\star }\mid
\lambda _{0}^{\star }\subset \eta \}$ we have:%
\begin{equation*}
\varphi ^{\star }(\lambda )=r^{\star }.
\end{equation*}

But $r=r^{\star }$ and $\lambda _{0}^{\star }\subset \Lambda $ by
construction. So, since $\Lambda \in \mathcal{F}^{\star }$, $\varphi ^{\star
}(\Lambda )=r.$

$\square $

\bigskip 

Having defined the $\Lambda $-limit, from now on we use the symbol $\ast $
to denote the extensions of objects in $\mathbb{U}$ in the sense of $\Lambda 
$-limit (not to be confused with the extensions obtained by applying the
star map $\star :$ e.g., the field $\mathbb{K=R}^{\ast }$ is a subfield of $%
\mathbb{R}^{\star }).$

We observe that, given a set $S$ in $\mathbb{U}$, its hyperfinite extension
(in the sense of the $\Lambda $-limit) is%
\begin{equation*}
S^{\circ }=\lim\limits_{\lambda \uparrow \mathbb{U}}(S\cap \lambda
)=S^{\star }\cap \Lambda
\end{equation*}%
and we use this observation to prove that, given a set of functions $%
V(\Omega )$, by posing%
\begin{equation*}
\widetilde{V}(\Omega )=Span(V(\Omega )^{\circ })=Span(V(\Omega )^{\star
}\cap \Lambda )
\end{equation*}%
we obtain the set of ultrafunctions generated by $V(\Omega )$.

The only nontrivial fact to prove is that, for every function $f\in V(\Omega
),$ its natural extension $f^{\ast }$ is an ultrafunction. First of all, we
observe that, by definition, $f^{\ast }=f^{\star }.$ Also, since $f\in
V(\Omega ),$ by transfer it follows that $f^{\star }\in $ $V(\Omega )^{\star
}$. And, by our choice of $\Lambda ,$ we also have that $f^{\star }\in
\Lambda $ since, by construction, $\{f\}^{\star }=\{f^{\star }\}\subset
\Lambda .$

It remains to prove the coherence of the axioms ($\Lambda $-1)\textsf{,}($%
\Lambda $-2),($\Lambda $-3) combined with $FTA$.

\begin{theorem}
The set of axioms $\{(\Lambda $-$1)$\textsf{,}$(\Lambda $-$2)$,$(\Lambda $-$%
3),FTA\}$ is consistent.
\end{theorem}

\textbf{Proof. }The basic idea is to chose an hyperfinite set $\Lambda \in
\bigcap_{\lambda \in \mathcal{F}}F_{\lambda }$,where $F_{\lambda }$ is
defined in Theorem \ref{catullo} (which automatically ensures the
satisfaction of ($\Lambda $-1)\textsf{,}($\Lambda $-2),($\Lambda $-3)), with
one more particular property that will ensure the satisfaction of $FTA.$

We start by considering a generic hyperfinite set $\Lambda ^{\prime }\in
\bigcap_{\lambda \in \mathcal{F}}F_{\lambda }$ and we let 
\begin{equation*}
B^{\prime }=\{e_{i}(x)|i\in I\}
\end{equation*}%
be any hyperfinite basis for $Span(V^{\mathfrak{F}}(\mathbb{R})^{\star }\cap
\Lambda ^{\prime }).$ Now we pose 
\begin{equation*}
B=\{\mathfrak{F}^{j}(e_{i}(x)):0\leq j\leq 3,i\in I\}\cup \{\overline{%
\mathfrak{F}^{j}(e_{i}(x))}:0\leq j\leq 3,i\in I\},
\end{equation*}

where $\mathfrak{F}$ denotes the Fourier transform. Since $\mathfrak{F}%
^{4}=id,$ we have that $B$ is closed by Fourier transform and complex
conjugate. We now pose%
\begin{equation*}
\Lambda =\Lambda ^{\prime }\cup B
\end{equation*}

and it is immediate to prove that, with this choice, $FTA$ is ensured,
because $B$ is a set of generators for $\widetilde{V^{\mathfrak{F}}(\mathbb{R%
})}$ closed by Fourier transform and complex conjugate.

$\square $

\end{document}